\documentclass[a4paper]{amsart}
\usepackage{amssymb,amsthm}
\address{Department of algebra,
      Faculty of Mechanics and Mathematics, 
      Moscow State University.
      Moscow 119991, Leninskie gory, MSU.
      \tt stas.mail@mail.ru}

\def\[{\left[}
\def\]{\right]}
\def\({\left(}
\def\){\right)}

\let\^\widehat
\let\~\widetilde

\let\=\overline
\let\ge\geqslant
\let\le\leqslant
\let\r\overrightarrow
\let\l\overleftarrow
\newtheorem*{lemma0}{Lemma 0}

\newtheorem{lemma}{Lemma}
\newtheorem*{theorem}{Theorem}

\author{Stanislav R. Safin}
\title{Powers of sets in free groups
}

\begin{document}

\vglue-12mm

\maketitle

\vskip-12mm

\noindent

\begin{abstract}
\noindent
We prove that $|A^n|\geqslant c_n\cdot|A|^{\[\frac{n+1}{2}\]}$ for any 
finite subset $A$ of a free group if $A$ contains at least two 
noncommuting elements, where $c_n>0$ are constants not depending on 
$A$.  Simple examples show that the order of these estimates are the best 
possible for each $n>0$.  
\end{abstract}

%%%%%%%%%%%%%%%%%%%%%%%%%%%%%%%%%%%%%%%%%%%%%%%%%%%%%%%%
\section{Introduction}

We use the following notation:
$F_m$ is the free group of rank $m$;
$|A|$ is the cardinality of a set $A$;
$AB$  is the set of all products of the form $ab$,
where $a\in A$ and $b\in B$;
$A^n$ is the set of all products of the form $a_1\dots a_n$,
where $a_i\in A$;
$|a|$ is the length of a 
word $a$ in some alphabet;
%% ������ $a$ ᢮������ ��㯯�;
$[x]$ is the largest integer not exceeding a real number $x$.

\smallskip

Mei-Chu Chang \cite{Chang} proved that there exist constants
$c,~\delta>0$ such that
$|A^3|>c\cdot|A|^{1+\delta}$
for any
finite subset~$A$ of the group ${\rm SL}_2(\mathbb C)$ 
not contained in any virtually abelian subgroup.
In particular, this estimate is valid for any finite subset
of a free group not contained in any cyclic subgroup.
A.~A.~Razborov \cite{Razborov}
improved Chang's estimate in this special case: 
there exists a constant $c>0$
such that $|A\cdot A\cdot A|\geqslant\frac{|A|^2}{(\log|A|)^c}$
for any finite subset $A$ of a
free group not contained in any cyclic subgroup.

It is easy to see that the squares of subsets satisfy no nontrivial
analogues of Chang and Razborov's inequality; the best 
possible estimate is linear in this case.

The following theorem shows, in particular, that the logarithm in 
Razborov's estimate can be removed.

\begin{theorem}
There exist constants $c_n>0$ such that, for any finite
subset $A$ of a free group not contained in any cyclic
subgroup, we have
$|A^n|\geqslant c_n\cdot|A|^{[\frac{n+1}{2}]}$
for all positive integers $n$.
\end{theorem}

It is easy to show that if $A=\{x,y,y^2,\dots,y^k\}$, then
$\left|A^n\right|=O\(k^{\[n+1\over2\]}\)$
for each fixed $n\ge1$,
where $x$ and $y$ are free generators of a free group.
Indeed, for $n=1$ and $n=2$, the assertion is obvious.
Further, by induction,
$
|A^n|=|A^{n-2}A^2|\le O\(k^{\[n-1\over2\]}\)\cdot|A^2|=
O\(k^{\[n-1\over2\]}\)\cdot O(k)=O\(k^{\[n+1\over2\]}\).
$
This simple example shows that
the order of our estimates is the best possible for each $n$ .

Note also that the proof below is significantly simpler than
the 
argument from \cite{Razborov}; but, on the other hand, in
\cite{Razborov}, a more general fact about the number of elements 
in the 
product of three (possibly, different) subsets of free group
was proven.

In commutative case, i.e., when the set $A$ is contained in
a
cyclic group, the situation is quite different. 
A survey of results on 
this subject can be found in \cite{Razborov} and in the 
literature cited therein.

\medskip
The author thanks A. A. Klyachko for setting the problem,
attention to the work, and several valuable
remarks.

%%%%%%%%%%%%%%%%%%%%%%%%%%%%%%%%%%%%%%%%%%%%%%%%%%%%%%%%
\section{Auxiliary lemmata}

The \emph{period},
or the \emph{left period}, of a word
$w$ in some alphabet is a
nonempty
%% , 横����᪨ ��᮪��⨬�� 
word $\~w$ not being a proper power such
that
the word $w$
has the form
$w=\~w^s\r w$, where $s\geqslant 2$ and
$\r w$ is a beginning of the word $\~w$
($\r w$ may be empty but must not coincide with $\~w$).
The word $\r w$ is called
the \emph{tail} (or the \emph{right tail}) of the word~$w$.
It is well known that the period and the tail of a periodic word 
are uniquely determined.

The \emph{right period}
of a word
$w$ in some alphabet is a nonempty
%%, 横����᪨ ��᮪��⨬�� 
word~$\^w$ not being a proper power
such
that
the word~$w$
has the form
$w=\l w\^w^s$, where $s\geqslant 2$ and
$\l w$ is an end of the word~$\^w$
($\l w$ may be empty but must not coincide with~$\^w$).
The word $\l w$ is called
the \emph{left tail} of the word~$w$.
It is well known that
the word has a right period if and only if it has
a left period (such words are called \emph{periodic}) and  lengths of periods are equal.
%%длины периодов совпадают.

In the following lemma, we collect some known properties of periodic
words.

\begin{lemma0}
If two periodic words have the same left periods and 
the same right
periods, then the tails of these words coincide too.
\newline
If, for some words, $u_1vw_1=u_2vw_2$
%%$=u_ivw_i$, 
%% ��᮪��⨬�, ������ 
and
$0<|u_2|-|u_1|\le{1\over2}|v|$
{\rm(}i.e., a word contains two occurrences of a word $v$ and these
occurrences intersect in a word of length at least
${1\over2}|v|${\rm)},
then the word $v$ is periodic and $u_2$ ends with $\~v$.
\end{lemma0}

\noindent
{\bf Proof.}
Let us prove the first assertion. Let's suppose that  $\~u=\~v$ and $\^u=\^v$. Without loss of generality, let's suppose $\r v=\r u t$, $\~u=\r utw$, where $\r utw$ - are reduced word. The word $u$ is ending with $tw\r u$, and the word $v$ is ending with $w\r ut$, but because these are the left and the right periods, it means they are the same. Because only  powers of one and the same word can commute \cite{Nielsen-Schreier-Th}, we see, that or $\~u$ is a proper power, or one of  them: $w\r u$ or $t$ is empty. $\~u$ can not be proper power because of the determination of the power. If $w\r u$ is empty, then $w$ is empty too, but it is impossible because of  determination of the tail  $\r v=\~v$. If $t$ is empty, then $\r u=\r v$, and we have proved our assertion.

\smallskip

The proof of the second assertion we leave to the reader as
an exercise.
In~\cite{Razborov}, this assertion (but with rougher estimate) 
was called the second overlapping lemma.

\begin{lemma}
% ����� 1.
For any finite set $A\subset F_2$, 
there exist a word $u\in F_2$ and sets 
${A_0,B_0\subseteq uAu^{-1}}$ 
such that $|A_0|,|B_0|\geqslant \frac1{62}|A|$ and

1) for any $a\in A_0$ and $b\in B_0$, the words $ab$ and $ba$ are reduced;

2) $|b|\ge|a|$
for all $b\in B_0$ and $a\in A_0$.
\end{lemma}

\noindent
\textbf{Proof.}
Let $x_1$ and $x_2$ be free generators of the group $F_2$, and let 
$e$ denote the empty word. We prove the lemma by induction on 
the sum of lengths of words 
from the set $A$.  
Let us decompose $A$ into the union of the 16 disjoint subsets 
$$ 
A_{(x,y)}=\{\hbox{words from $A$ beginning with $x$
and ending with $y$}\},
\quad\hbox{where $x,y\in\{x_1,x_1^{-1},x_2,x_2^{-1}\}$}.
$$
(If $A$ contains the empty word, then we include it in
$A_{(x_1,x_1)}$.)

\smallskip
\noindent
{\bf Case 1.}
$|A_{(x,y)}|\ge\frac{1}{31}|A|$
for some non mutually inverse $x,y$.
In this case, we 
put $A_0$ to be
the set of
$\[{|A_{(x,y)}|+1\over2}\]$ shortest words from
$A_{(x,y)}$
and put $B_0$ to be the set of $\[{|A_{(x,y)}|+1\over2}\]$
longest words from $A_{(x,y)}$.
Clearly, these sets $A_0$
and $B_0$ are as required (with $u=e$).

\smallskip
\noindent
{\bf Case 2.}
$|A_{(x,x^{-1})}|\ge\frac{1}{31}|A|\le|A_{(y,y^{-1})}|$
for some different $x,y$.
Without loss of generality, we assume that a mean 
(by length)
word from $A_{(x,x^{-1})}$
is not longer than a mean word from $A_{(y,y^{-1})}$.
In this case, we put
$A_0$ to be the set of $\[{|A_{(x,x^{-1})}|+1\over2}\]$ shortest words 
from $A_{(x,x^{-1})}$, and put $B_0$ to be the set of 
$\[{|A_{(y,y^{-1})}|+1\over2}\]$ longest words from $A_{(y,y^{-1})}$.  
Clearly, these sets $A_0$
and $B_0$ are as required (with $u=e$).

\smallskip
\noindent
{If the conditions of neither case 1 nor case 2 hold,}
then, obviously, for some letters $x$,
$$
|A_{(x,x^{-1})}|>\(1-\frac{15}{31}\)|A|=\frac{16}{31}|A|>{1\over2}|A|
$$
and, therefore, the total length of words
of the 
set $x^{-1}Ax$ is less than the total length of words
of the 
set $A$. The application of the induction hypothesis completes proof.

\begin{lemma}
% ����� 2.
Suppose that $U,V,W\subset F_2$ and all products $UVW$ are reduced.
If $|v|\ge|u|$ for all $u\in U$ and $v\in V$,
then either $|UVW|\geqslant \frac16|U|\cdot |W|$ or
all words from $V$ are periodic with the same
period.
{\rm Similarly,}
if $|v|\ge|w|$ for all $w\in W$ and $v\in V$,
then either $|UVW|\geqslant \frac16|U|\cdot |W|$ or
all words from $V$ are periodic with the same
right period.
\end{lemma}

This lemma
is proven in the last section. Now, 
let us derive the theorem from these two lemmata.

%%%%%%%%%%%%%%%%%%%%%%%%%%%%%%%%%%%%%%%%%%%%%%%%%%%%%%%%
\section{Proof of the theorem}

The words from $A$ have only finitely many different letters.
Therefore, $A$
is contained in a free group $F_m$ of finite rank.
Since, as is known,
$F_m$ embeds into $F_2$
\cite{Nielsen-Schreier-Th}, we can assume that $m=2$.
Clearly, it is sufficient to prove the assertion of the theorem 
for odd $n$; so, we assume that $n=2k+1$.

Applying Lemma 1 to $A$, we obtain
sets $A_0, B_0\subseteq u^{-1}Au$. Note
that $|(uAu^{-1})^{n}|=|A^{n}|$; therefore, we can assume that $u=e$.

First, consider the case where the set $B_0$
consists of periodic words with the same
period $p$ and the same tail $t$:
$$
B_0=\{p^{n_1}t, p^{n_2}t,\dots\}.
\eqno(*)
$$

If $t=e$, then there exists a word $b\in A$ not commuting with
$p$, because by condition,
the set $A$ contains noncommuting elements.
The group generated by $p$ and $b$ is generated by them freely,
because any two noncommuting words freely
generate a free group of rank 2 (see \cite{Nielsen-Schreier-Th}).
Therefore, all products of the 
form $u_1b\dots u_kbu_{k+1}$, where $u_i\in B_0$, are different and, hence,
$
|A^n|\ge |(B_0b)^kB_0| = |B_0|^{k+1}\geqslant O(|A|^{k+1}).
$

If $t\ne e$, then we estimate $|B_0^n|$. Since
$p$ and $t$ do not commute
(because, in a free group, only powers of the same element commute
\cite{Nielsen-Schreier-Th},
the period $p$ is not a proper power, and $|t|<|p|$),
they freely generate a 
free group of rank~2.
Therefore, all words $\prod\limits_{i=1}^{n}u_i$,
where $u_i\in B_0$, are different and $|B_0^n|\geqslant O(|A|^n)$.

Now, consider the case where
$B_0$ is not of the form $(*)$.
Let us 
prove the inequality
$$
{|(A_0B_0)^kA_0|\geqslant O(|A|^{k+1})}
$$
by induction on
$k$. For $k=0$, the assertion is obvious.

By Lemma 0, either not all left periods of words from $B_0$ coincide, 
or not all right periods of words from $B_0$ coincide. Without loss of 
generality, we assume that the left periods do not coincide.

Applying Lemma 2 to the sets $U=A_0$, $V=B_0$, and
$W=(A_0B_0)^{k-1}A_0$, we
obtain the required inequality
$$
|A_0\cdot B_0\cdot(A_0B_0)^{k-1}A_0|
\geqslant O(|A_0|\cdot |(A_0B_0)^{k-1}A_0|)\geqslant
O(|A|^{k+1}).
$$

%%%%%%%%%%%%%%%%%%%%%%%%%%%%%%%%%%%%%%%%%%%%%%%%%%%%%%%%
\section{Proof of Lemma 2}

\begin{lemma}
% ����� 3.
Suppose that words 
$u_1,u_2,u_3,v,w_1,w_2,w_3\in F_2$ are such that 
$u_1vw_1=u_2vw_2=u_3vw_3$,
all words $u_ivw_i$ are reduced,
$u_i$
are 
pairwise different,  and $|v|\geqslant |u_i|$. 
Then the word $v$ is periodic with 
period $\~v$ and one of the words $u_i$ ends with~$\~v$.  
\end{lemma}

\noindent
\textbf{Proof.}
We see that the word $f=u_1vw_1=u_2vw_2=u_3vw_3$
has three
occurrences of the subword $v$.
Without loss of generality, we can assume that $|u_1|<|u_2|<|u_3|$.
Since $|v|\ge|u_i|$,
any two of these three occurrences of 
$v$ either intersect or, at least, are adjacent to 
each
other (i.e. the $k$-th letter of $f$ is the end of one occurrence of $v$  
and
the $(k+1)$-th letter of $f$ is the beginning of the other occurrence of 
$v$).  Therefore, the second occurrence of~$v$ is completely covered by 
the first and third occurrences of~$v$ and, hence, the second occurrence 
of~$v$ intersects with one of the other occurrences (say, the first) in a 
subword of length at least ${1\over2}|v|$.  By the second assertion of 
Lemma 0, this means that the word $v$ is periodic and the word $u_2$ ends 
with~$\~v$.

\begin{lemma}
% ����� 4.
Suppose that sets $U$, $V$, and $W$ satisfy conditions of Lemma 2 and
$v\in V$.
If
$U$ contains no words ending with the period of~$v$
or the word $v$ is nonperiodic,
then
$|UvW|\geqslant \frac12|U|\cdot |W|$.
\end{lemma}

\noindent
\textbf{Proof.}
Lemma 3 shows that no word has more than two representations as 
a product $uvw$,
where $u\in U$ and $w\in W$.
Therefore,
$|UvW|\geqslant \frac12|U|\cdot|W|$, as required.

\begin{lemma}
% ����� 5.
Suppose that $\alpha$ and $\beta$ are nonempty cyclically reduced 
words that are not proper powers, $|\beta|>|\alpha|$, \ 
$a,~b\in F_2$ are words ending with $\beta^2$, and $a$ ends with 
$\alpha^{s}$. Then $b$ ends with $\alpha^{s}$ too.

\end{lemma}

\noindent
\textbf{Proof.}
If
$|\alpha^s|\leqslant|\beta^2|$,
then
$a$ and $b$ end with $\beta^2$; so, they end with $\alpha^{s}$,
as required.
If $|\alpha^{s}|>|\beta^2|$,
then the end $\beta^2$ of~$a$
has two different right periods,
$\alpha$ and $\beta$;
this
contradiction completes the proof.

\begin{lemma}
% ����� 6.
Suppose that sets $U$, $V$, and $W$ satisfy conditions of Lemma 2 and
there exists a periodic word~$v\in V$ with
period $\~v$ such that all words in $U$ end with
$\~v^q$ and
none of them ends with $\~v^{q+1}$, where $q\geqslant 1$. Then
$|UvW|\geqslant \frac12|U|\cdot|W|$.
\end{lemma}

\noindent
\textbf{Proof.}
This lemma follows immediately from Lemma 4 in which for 
$U$, $V$, and~$W$  the sets
$U\~v^{-q}$, $\~v^qV$, and~$W$, respectively, 
are taken
and
the word $\~v^qv$ is taken for $v$.

\medskip
\noindent
\textbf{Proof of Lemma 2.}
Clearly, it is sufficient to prove the first assertion.
By Lemma 4, we may assume that the set $V$ consists of periodic words. 

Suppose that $V$ has two words $v_1$ and $v_2$ with different
periods $\~v_1$ and $\~v_2$.
Suppose that $|\~v_1|\ge|\~v_2|$.
Applying Lemma~4 once again, we obtain that either
the required inequality holds or
there exists a set $U_0\subseteq U$ such that
$|U_0|\geqslant\frac{1}{3}|U|$
and all words in $U_0$
end with both $\~v_1$ and $\~v_2$.
This means, in particular, that $\~v_1$
ends with $\~v_2$ and $|\~v_1|>|\~v_2|$
(because $\~v_1\ne\~v_2$).

Let $U_{00}$ be the set of words from $U_0$ ending with
$\~v_2^2$. By Lemma 5, all words from
$U_{00}$ end with
$\~v_1^s$, and none of them ends with $\~v_1^{s+1}$.
By Lemma 6 (in which the role of $U$ is played by $U_{00}$ and the role
of $v$ is played by $v_1$),
we have
$$
|U_{00}v_1W|\ge \frac12|U_{00}|\cdot|W|.
\eqno(**)
$$

On the other hand,
in the set
$U_0\setminus U_{00}$, all words end with
$\~v_2$, but none ends with $\~v_2^2$.
Therefore, by Lemma 6 (in which the role of $U$ is played by
$U_0\setminus U_{00}$, and the role of~$v$ is played by~$v_2$), we
have
$$
|(U_0\setminus U_{00})v_2W|\ge \frac12|U_0\setminus U_{00}|\cdot|W|.
\eqno({**}*)
$$
Inequalities $(**)$ and $({**}*)$ imply 
$|U_0\cdot(\{v_1\}\cup \{v_2\})\cdot W|\ge \frac12|U_0|\cdot|W|=
\frac16|U|\cdot|W|$, as required.

%%%%%%%%%%%%%%%%%%%%%%%%%%%%%%%%%%%%%%%%%%%%%%%%%%%%%%%%%%

\end{document}